\RequirePackage{fix-cm}
  
\documentclass[smallextended]{svjour3}       

\smartqed  

\usepackage{graphicx}

\usepackage[english]{babel} 
\usepackage{amssymb}
\usepackage{amsmath}
\usepackage{txfonts}
\usepackage{mathdots}
\usepackage[classicReIm]{kpfonts}
\usepackage{graphicx}
\usepackage{amscd}
\usepackage{comment}
\usepackage[margin=1.3in]{geometry}
\allowdisplaybreaks
\journalname{}

\begin{document}

\title{Riemann's Last Theorem 

}

\author{Aric Behzad-Canaanie}

\institute{Aric Behzad-Canaanie \at
	Aecom/Nasa Ames Research Center \\
    \emph{M/S 213-8 Bldg. 213, Moffett Field, CA 94035} \\ 
	American Mensa \\
	\emph{1315 Brookside Dr., Hurst, TX 76053\\
	\email{behzad.canaanie@member.mensa.org}\\   
	}
	} 

\date{\today}

\maketitle

\section*{Abstract}
\label{abst}
	
\noindent The central idea of this article is to introduce and prove a special form of the zeta function as proof of Riemann’s last theorem. The newly proposed zeta function contains two sub functions, namely $f_{\mathrm{1}}\left(b,s\right)$ and $\mathrm{\ }f_{\mathrm{2}}\left(b,s\right)$. The unique property of $\mathrm{\zetaup }\left({s}\right)=f_{\mathrm{1}}\left(b,s\right)\mathrm{-}f_{\mathrm{2}}\left(b,s\right)$ is that as $b$ tends toward $\mathrm{\infty },$ the equality $\mathrm{\zetaup }\left(s\right)\mathrm{=}\mathrm{\ }\mathrm{\zetaup }\left(1-s\right)$ is transformed into an exponential expression  for the zeros of the zeta function. At the limiting point, we simply deduce that the exponential equality is satisfied if and only if $\mathfrak{R}\left(s\right)\mathrm{=}\frac{\mathrm{1}}{\mathrm{2}}$. Consequently, we conclude that the zeta function cannot be zero if $\mathfrak{R}\left(s\right)\mathrm{\neq }\frac{\mathrm{1}}{\mathrm{2}}$, hence proving Riemann's last theorem.\\

%\end{abstract}

\section*{Introduction}
\label{intro}
 \noindent Riemann's last theorem (RSLT) is very well known and providing a traditional introduction is not the most efficient way to achieve the objective of this work. This subject is well explained in countless publications [\textit{(\Ref1), (\Ref2), (\Ref13), (\Ref14)}] and the assumption is that the reader understands this topic; therefore, instead of explaining RSLT, please consider a brief history of a singular number known as zero. This number reformed our understanding of mathematics, and it’s responsible for many advancements and achievements in mathematics. For example, the solution for an equation or function is the zero(s) of that equation or function. It’s not so hard to see that the entire calculus is the solution to  $\frac{0}{0}$.  It’s mysterious why this number reveals so much information. For instance, by knowing the zeros of any system, we can construct an accurate model of that system. This work is based on the behavior of functions close to infinity that has been rewritten in a more known and familiar methods. And the paragraph below has been kept as a testimony to the greatness of Riemann’s work. Where he revealed an important connection between the zeros of the zeta function and the prime numbers.\\ 

\noindent The English philosopher Adelard of Bath is known to have adopted the numeral zero to the number system in the 12th century [\textit{(5)}]. This critical introduction has initiated a journey for this number with respect to the roots of mathematical functions. Over the centuries, the number zero has become the most important source of challenges in mathematics.\textit{ }One such challenge involves proving a conjecture, which states that the real part of the nontrivial zeros of the zeta function [\textit{(1), Chap. X}] is 1/2. This conjecture was first published by Georg Friedrich Bernhard Riemann in 1859 [\textit{(2), Chap. I}] and is important because it provides a solid and irreplaceable foundation to gain information about the distribution of prime numbers [\textit{(4)}]. In this work, we present proof of this conjecture.\\

\section*{Assumptions}
\label{Assum}
	 $\mathrm{1}\ \mathrm{\le }\mathrm{\ }{x}\ \mathrm{\in }\ \mathbb{R}$$\mathrm{;}$ ${n= }\left\lfloor x\right\rfloor\mathrm{\in }\mathrm{\ }$$\mathrm{\mathbb{N}}$$\mathrm{;}$$\mathrm{\ }\mathrm{\Gamma }\left(s\right)\mathrm{\neq }\ \mathrm{0}$,\ $ s\ \mathrm{\in }\mathrm{\ }$$\mathrm{\mathbb{C}};\ \left\{\ \right\}\neq\left\{ 0 \right\}\neq\left\{0,0 \right\}\neq \dots$.\\
	
\section*{Boundary}
\label{Bound}	
The functions, expressions, and Equations [\eqref{2}--\eqref{6}] are restricted to the critical region $\left(0<\Re\left(s\right)<1\right)$.\\

\section*{Consideration}
\label{Consid}
The nontrivial roots of the zeta function are considered for validating all relevant statements and conclusions.\\

\section*{Definitions}
\label{Defin}
For $\mathrm{s=}\mathrm{\sigmaup }\mathrm{+it}\mathrm{\ }$$\mathrm{;}$$\mathrm{\ }\mathrm{\sigmaup }\mathrm{,\ t}\ \mathrm{\in }\ \mathbb{R}\mathrm{;}\mathrm{\ }\mathrm{i=}\sqrt{\mathrm{-1}}$, the zeta function [\textit{(1), Chap. I}] is defined as 
\begin{equation}\label{1}
	\zeta\left(s\right)=\sum_{n=1}^{\infty}\frac{1}{n^s}=\sum_{n=1}^{\infty}\frac{1}{n^{\sigma+it}}
\end{equation}

\noindent This function converges for $\mathfrak{R}\left(\mathrm{s}\right)\mathrm{>1}$ and meromorphically continues over the entire complex plane with a simple pole residue of 1 at $\mathrm{s=1}$.\\

\noindent The zeta function then satisfies the functional equation shown below [\textit{(1), Chap. II}]:
\begin{equation}\label{1-1}
	\zeta\left(s\right)=\frac{2^s}{\pi^{1-s}}sin{\left(\ \frac{\pi}{2}s\right)}\Gamma\left(1-s\right)\zeta\left(1-s\right)
\end{equation}

\noindent Equation \eqref{1} can now be rewritten [\textit{(1), Chap. II}] for the critical region $\mathrm{(0< \Re\left(s\right)=\sigmaup<1)}$ as

\begin{equation}\label{1-2}
	\zeta\left(s\right)=\frac{1}{s-1}+\sum_{n=1}^{\infty}\int_{n}^{n+1}\left(\frac{1}{n^s}-\frac{1}{x^s}\right)dx=s\left(\frac{1}{s-1}-\int_{1}^{\infty}\frac{x-\left\lfloor x\right\rfloor}{x^{s+1}}dx\right)
\end{equation}\\
\section*{Riemann's last theorem}
\label{RSLT}

The real part of the nontrivial zeros $\left(\left\{\ s\ |\  0<\mathfrak{R}\left(s\right)<1,\ \zeta \left(s\right)=0\ \right\}\right)$  of the zeta function is 1/2.\\
\section*{Proof}
\label{proof}

\noindent By applying the summation and integration properties to the left-hand-side function of \eqref{1-2}, we can say that
\begin{align*}
	\zeta\left(s\right)&=\frac{1}{s-1}+\sum_{n=1}^{\infty}\int_{n}^{n+1}\left(\frac{1}{n^s}-\frac{1}{x^s}\right)dx\\
	&=\frac{1}{s-1}+\sum_{n=1}^{\infty}\int_{n}^{n+1}\left(\frac{1}{n^s}\right)dx-\sum_{n=1}^{\infty}\int_{n}^{n+1}\left(\frac{1}{x^s}\right)dx\\
	&=\frac{-1}{1-s}+{\left(\sum^{\infty}_{n\mathrm{=1}}{\left(\frac{\mathrm{1}}{n^s}\right)}\ -\int^{\infty}_1{\frac{1}{x^s}dx}\right)}\\
	&=\frac{-1}{1-s}+{\mathop{\mathrm{lim}}_{b\mathrm{\to }\mathrm{\infty }} \left(\sum^b_{n\mathrm{=1}}{\left(\frac{\mathrm{1}}{n^s}\right)}-\ \frac{b^{1-s}}{1-s}+\frac{1}{1-s} \right)}		
\end{align*}
\noindent Therefore, we obtain	
\begin{equation}\label{2}
	\mathrm{\zetaup }\left(s\right)={\mathop{\mathrm{lim}}_{b\mathrm{\to }\mathrm{\infty }} \left(\sum^b_{n\mathrm{=1}}{\left(\frac{\mathrm{1}}{n^s}\right)}\ -\ \frac{b^{1-s}}{1-s}\right)\ }
\end{equation}

\noindent Further, we can generalize the function on the far right of \eqref{1-2} as follows: 
\begin{align*}%\label{2-1}
	\zeta \left(s\right)&=\frac{s}{s-1}-s\int_{1}^{\infty}\frac{x-\left\lfloor x\right\rfloor}{x^{s+1}}dx\\
	&=\frac{s}{s-1}-\sum^{\infty}_{n\mathrm{=1}}{\left(s\int^{n+1}_n{\frac{x-n}{x^{s+1}}}dx\right)}\\
	&=\frac{s}{s-1}-\sum^{b-1}_{n\mathrm{=1}}{\left(s\int^{n+1}_n{\frac{x-n}{x^{s+1}}}dx\right)}-\underbrace{s\int^{\infty }_b{\frac{x-\left\lfloor x\right\rfloor }{x^{s+1}}}dx}_{O_1}\\
	&=\frac{s}{s-1}-s\int^b_1{\frac{x}{x^{s+1}}}dx-s\sum^{b-1}_{n\mathrm{=1}}{\left(\int^{n+1}_n{\frac{-n}{x^{s+1}}}dx\right)}-{O_1}\\
	&=\frac{-s}{1-s}+s\int^1_b{\frac{1}{x^s}}dx-s\sum^{b-1}_{n\mathrm{=1}}{n\left(\int^n_{n+1}{\frac{1}{x^{s+1}}}dx\right)}-{O_1}\\
	&=\frac{-s}{1-s}+\frac{s}{1-s}-\frac{s}{1-s}b^{1-s}+\sum^{b-1}_{n\mathrm{=1}}{\left(\frac{n}{n^s}-\frac{n}{{\left(n+1\right)}^s}\right)}-{O_1}
\end{align*}

\noindent Considering $\frac{-s\ }{1-s}b^{1-s}=\frac{b}{b^s}-\frac{b^{1-s}}{1-s}$, we have
\begin{align*}%\label{2-1-1}
	\zeta \left(s\right)&=\frac{b}{b^s}-\frac{b^{1-s}}{1-s}+\left(\frac{1}{1^s}-\frac{1}{2^s}+\frac{2}{2^s}-\frac{2}{3^s}+\dots +\frac{b-1}{{\left(b-1\right)}^s}-\frac{b-1}{b^s}\right)-{O_1}\\
	&=\frac{b}{b^s}-\frac{b^{1-s}}{1-s}+\sum^{b-1}_{n\mathrm{=1}}{\left(\frac{\mathrm{1}}{n^s}\right)}-\frac{b}{b^s}+\frac{1}{b^s}-\underbrace{s\int^{\infty }_b{\frac{x-\left\lfloor x\right\rfloor }{x^{s+1}}}dx}_{O_1}
\end{align*}
\noindent Thus, we obtain a general (valid for b = 1, 2, 3...) form of the right-hand side (RHS) function of \eqref{1-2} as
\begin{align}\label{2-1-2}
	\zeta \left(s\right)=\sum^b_{n\mathrm{=1}}{\left(\frac{\mathrm{1}}{n^s}\right)}\ -\frac{b^{1-s}}{1-s}-s\int^{\infty }_b{\frac{x-\left\lfloor x\right\rfloor }{x^{s+1}}}dx\ ,\ b\ \mathrm{\in }\mathrm{\ }\mathbb{N}
\end{align}

\noindent Note that $\mathrm{b=1}$ provides the RHS function of \eqref{1-2}. By considering the limit of \eqref{2-1-2} as $b$ tends toward infinity, \eqref{2} can be proved by an alternative method.\\

\noindent We can express \eqref{1-2} for $\mathrm{\zetaup }\left(1-s\right)$, where$\ \mathrm{s=}\mathrm{\sigmaup }\mathrm{+it}$, as 
\begin{align}\label{2-2}
	\begin{split}
		\zeta \left(1-s\right)&=\frac{1}{-s}+\sum^{\infty }_{n=1}{\int^{n+1}_n{\left(\frac{1}{n^{1-s}}-\frac{1}{x^{1-s}}\right)}}dx\\
		&=\left(1-s\right)\left(\frac{1}{-s}-\int^{\infty }_1{\frac{x-\left\lfloor x\right\rfloor }{x^{2-s}}}dx\right)
	\end{split}
\end{align}

\noindent Therefore, considering \eqref{2},\ \eqref{2-1-2}, and \eqref{2-2}, we can derive that
\begin{align}\label{2-3}
	\begin{split}
		{\mathop{\mathrm{lim}}_{b\to \infty } \left(\sum^b_{n\mathrm{=1}}{\left(\frac{\mathrm{1}}{n^s}\right)}\ -\frac{b^{1-s}}{1-s}\right)\ } &=\sum^b_{n\mathrm{=1}}{\left(\frac{\mathrm{1}}{n^s}\right)}\ -\frac{b^{1-s}}{1-s}-s\int^{\infty }_b{\frac{x-\left\lfloor x\right\rfloor }{x^{s+1}}}dx \ ,\ b\ \mathrm{\in }\mathrm{\ }\mathbb{N}\\ 
		{\mathop{\mathrm{lim}}_{b\to \infty } \left(\sum^b_{n\mathrm{=1}}{\left(\frac{\mathrm{1}}{n^{1-s}}\right)}-\frac{b^s}{s}\right)\ }&=\sum^b_{n\mathrm{=1}}{\left(\frac{\mathrm{1}}{n^{1-s}}\right)}\ -\frac{b^s}{s}-\left(1-s\right)\int^{\infty }_b{\frac{x-\left\lfloor x\right\rfloor }{x^{2-s}}}dx\ ,\ b\ \mathrm{\in }\mathrm{\ }\mathbb{N}
	\end{split}
\end{align}
\noindent Considering \eqref{1-1} and the symmetry of the zeta function, we can say that if$\ \mathrm{\zetaup }\left(s\right)=0$, then $\ \overline{\mathrm{\zetaup }\left(1-s\right)}\ \ $will also be equal to zero. Therefore, we may assume that the roots of the zeta functions would occur at $\mathrm{\zetaup }\left(s\right)=\overline{\mathrm{\zetaup }\left(1-s\right)}$. Hence, we have 
\begin{align}\label{2-4}
	\begin{split}
		\mathrm{\zetaup }\left(s\right)=0\ \ & \mathrm{\Rightarrow }\ \mathrm{\zetaup }\left(s\right)=\overline{\mathrm{\zetaup }\left(1-s\right)}\mathrm{\ }
		\\ &\mathrm{\Rightarrow }\ {\mathop{\mathrm{lim}}_{b\mathrm{\to }\mathrm{\infty }} \left(\sum^b_{n=1}{\left(\frac{1}{n^s}\right)}-\frac{b^{1-s}}{1-s}\right)\ }={\mathop{\mathrm{lim}}_{b\mathrm{\to }\mathrm{\infty }} \left(\sum^b_{n\mathrm{=1}}{\left(\overline{\frac{1}{n^{1-s}}}\right)}-\overline{\left(\frac{b^s}{s}\right)}\right)\ } 
	\end{split}
\end{align}
\noindent Alternately, by considering the RHS functions of \eqref{2-3}, we have
\begin{align}\label{2-4-1}
	\begin{split}
		\mathrm{\zetaup }\left(s\right)=\overline{\mathrm{\zetaup }\left(1-s\right)}\mathrm{\ }\ \mathrm{\Rightarrow }
		&\sum^b_{n\mathrm{=1}}{\left(\frac{\mathrm{1}}{n^s}\right)} -\frac{b^{1-s}}{1-s}-{s}\int^{\mathrm{\infty }}_b{\frac{x-\left\lfloor x\right\rfloor }{x^{s+1}}{dx}}\\
		=&\sum^b_{n\mathrm{=1}}{\left(\overline{\frac{1}{n^{1-s}}}\right)}-\overline{\left(\frac{b^s}{s}\right)}-\overline{\left(1-s\right)\int^{\mathrm{\infty }}_b{\frac{x-\left\lfloor x\right\rfloor }{x^{2-s}}dx}}
	\end{split}
\end{align}

\noindent Note that \eqref{2-4-1} is valid for b = 1, 2, 3, {\dots}. Thus, we note that the limit of \eqref{2-4-1} as \textit{b} approaches infinity asserts \eqref{2-4} by an alternative method.\\

\noindent Moving the subfunctions $\sum^b_{n\mathrm{=1}}{\left(\overline{\frac{1}{n^{1-s}}}\right)}$ to the RHS and $-\frac{b^{1-s}}{1-s}\ $to the LHS gives 
\begin{align}\label{2-4-2}
	\begin{split}
		\mathrm{\zetaup }\left(s\right)=\overline{\mathrm{\zetaup }\left(1-s\right)}\mathrm{\ }\ \mathrm{\Rightarrow }&
		\sum^b_{n\mathrm{=1}}{\left(\frac{\mathrm{1}}{n^s}\right)}\ -\sum^b_{n\mathrm{=1}}{\left(\overline{\frac{1}{n^{1-s}}}\right)}-{s}\int^{\mathrm{\infty }}_b{\frac{x-\left\lfloor x\right\rfloor }{x^{s+1}} {dx}}\\
		=&\frac{b^{1-s}}{1-s}-\overline{\left(\frac{b^s}{s}\right)}-\overline{\left(1-s\right)\int^{\mathrm{\infty }}_b{\frac{x-\left\lfloor x\right\rfloor }{x^{2-s}}dx}}
	\end{split}
\end{align}

\noindent Then, taking the limit of \eqref{2-4-2} as $\mathrm{b}\mathrm{\to }\mathrm{\infty }$ gives
\begin{align}\label{2-4-3}
	\begin{split}
		\mathrm{\zetaup }\left(s\right)=0\ \ & \mathrm{\Rightarrow }\ \mathrm{\zetaup }\left(s\right)=\overline{\mathrm{\zetaup }\left(1-s\right)}\mathrm{\ }\ \\ & \mathrm{\Rightarrow }\ {\mathop{\mathrm{lim}}_{b\mathrm{\to }\mathrm{\infty }} \sum^b_{n=1}{\left(\frac{1}{n^s}-\left(\overline{\frac{1}{n^{1-s}}}\right)\right)}\ }={\mathop{\mathrm{lim}}_{b\mathrm{\to }\mathrm{\infty }} \left(\ \frac{b^{1-s}}{1-s}-\overline{\left(\frac{b^s}{s}\right)}\right)\ }
	\end{split}
\end{align}

%__________________________________________________________________________________________________________________________________________________*************************************

\noindent The simple and robust (addition, subtraction, multiplication, and division) steps below demonstrate a general technique to manipulate all infinite series of the form $1+\frac{1}{2^s}+\frac{1}{3^s}+{\dots}$ to produce finite values.

\begin{align}\label{2-5-e}
	\begin{split}
		1+\frac{1}{2^s}+\frac{1}{3^s}+{\dots} &=\sum^{\infty }_{n\mathrm{=1}}{\frac{\mathrm{1}}{n^s}}\\
		&=\sum^{\infty }_{n\mathrm{=1}}{\frac{\mathrm{1}}{n^s}}\left(1-\frac{2}{2^s}\right){\left(1-\frac{2}{2^s}\right)}^{-1}\\
		&=\left(\sum^{\infty}_{n\mathrm{=1}}{\left(\frac{\mathrm{1}}{n^s}\right)\left(1-\frac{2}{2^s}\right)}\right){\left(1-\frac{2}{2^s}\right)}^{-1}\\
		&=\left(\sum^{\infty}_{n\mathrm{=1}}{\left(\frac{\mathrm{1}}{n^s}\right)}-\frac{2}{2^s}\sum^{\infty}_{n\mathrm{=1}}{\left(\frac{1}{n^s}\right)}\right){\left(1-\frac{2}{2^s}\right)}^{-1}\\
		&=\left(\sum^{\infty}_{n\mathrm{=1}}{\left(\frac{\mathrm{1}}{n^s}\right)}-\sum^{\infty}_{n\mathrm{=1}}{\left(\frac{1}{{\left(2n\right)}^s}+\frac{1}{{\left(2n\right)}^s}\right)}\right){\left(1-\frac{2}{2^s}\right)}^{-1}\\
		&=\left(\sum^{\mathrm{\infty}}_{n\mathrm{=1}}{\left(\frac{1}{{\left(2n-1\right)}^s}+\frac{1}{{\left(2n\right)}^s}-\frac{1}{{\left(2n\right)}^s}-\frac{1}{{\left(2n\right)}^s}\right)}\right){\left(1-\frac{2}{2^s}\right)}^{-1}\\
		&=\left(\sum^{\mathrm{\infty}}_{n\mathrm{=1}}{\left(\frac{1}{{\left(2n-1\right)}^s}-\frac{1}{{\left(2n\right)}^s}+\frac{1}{{\left(2n\right)}^s}-\frac{1}{{\left(2n\right)}^s}\right)}\right){\left(1-\frac{2}{2^s}\right)}^{-1}\\
		&=\left(\sum^{\mathrm{\infty }}_{n\mathrm{=1}}{\left(\frac{1}{{\left(2n-1\right)}^s}-\frac{1}{{\left(2n\right)}^s}\ +0\right)}\right){\left(1-\frac{2}{2^s}\right)}^{-1}\\
		&=\left(\sum^{\mathrm{\infty }}_{n\mathrm{=1}}{\left(\frac{1}{{\left(2n-1\right)}^s}-\frac{1}{{\left(2n\right)}^s}\right)}\right){\left(1-\frac{2}{2^s}\right)}^{-1}\\
		&=\left(\sum^{\infty }_{n\mathrm{=1}}{\frac{{\mathrm{-}\mathrm{1}}^{n-1}}{n^s}}\right){\left(1-\frac{2}{2^s}\right)}^{-1}\\
		&=\frac{1}{1-2^{1-s}}\sum^{\infty }_{n\mathrm{=1}}{\frac{{\mathrm{-}\mathrm{1}}^{n-1}}{n^s}}\\
		&=\zeta \left(s\right)
	\end{split}
\end{align}

\noindent The term $\frac{1}{1-2^{1-s}}$ is regular for all $\mathrm{s}\mathrm{\ }\mathrm{\in }\mathrm{\ }\mathrm{\mathbb{C-}}\left\{1\right\}$, and the function $\sum^{\mathrm{\infty }}_{n\mathrm{=1}}{\frac{\mathrm{-}{\mathrm{1}}^{n-1}}{n^s}}$ is commonly known as the Dirichlet eta function that converges for all complex numbers with $\mathfrak{R}\left(s\right)>0$. Therefore, $\frac{1}{1-2^{1-s}}\sum^{\infty }_{n\mathrm{=1}}{\frac{{\mathrm{-1}}^{n-1}}{n^s}}$ [\textit{(1), Chap. II}] converges within the critical region $\mathrm{(0<\sigmaup<1)}$.\\ 

\noindent Consider replacing  $\sum^{\infty }_{n\mathrm{=1}}{\frac{\mathrm{1}}{n^s}}$ with $f_1\left(s\right),\ $ $\sum^{\infty }_{n\mathrm{=1}}{\frac{\mathrm{1}}{n^s}}\left(1-\frac{2}{2^s}\right){\left(1-\frac{2}{2^s}\right)}^{-1}$ with  $f_2\left(s\right)$, and so on; thus, we obtain the simplified version of \eqref{2-5-e} as follows: 
\begin{align*}%\label{2-5-e3}
	\begin{split}
		1+\frac{1}{2^s}+\frac{1}{3^s}+...=f_1\left(s\right)=f_2\left(s\right)=...=f_{11}\left(s\right)=\zeta \left(s\right)\ \
	\end{split}
\end{align*} 
\noindent Here, we can easily observe that \eqref{2-5-e} proves a one-to-one correspondence (mapping) or bijection between the series  $1+\frac{1}{2^s}+\frac{1}{3^s}+\textit{{\dots}}$ and $\zeta \left(s\right)$ in the critical region simply because for each s (input), all functions will generate one and only one series (output). In other words, all functions (function definition [\textit{(3), Chap. 5.1-page 66}]) in \eqref{2-5-e} are in bijection with s (input) and each other because none of the functions generate two or more series for each s. The main idea of \eqref{2-5-e} is to prove the one-to-one correspondence (mapping) or bijection, and later (it is difficult to show all these expressions simultaneously) in \eqref{2-5}, \eqref{2-5-1}, and \eqref{2-5-2}, we will see how to derive the simplified version of \eqref{2-5-e} in a manner such that all functions become regular for $\mathrm{\zetaup }\left(s\right)=\overline{\mathrm{\zetaup }\left(1-s\right)}$.\\

\noindent The significance of \eqref{2-5-e} is that it proves the bijection between the $\mathrm{\zetaup }\left(s\right)=\frac{1}{1-2^{1-s}}\sum^{\mathrm{\infty }}_{n\mathrm{=1}}{\frac{\mathrm{-}{\mathrm{1}}^{n-1}}{n^s}}$  and the series $\sum^{\mathrm{\infty }}_{n\mathrm{=1}}{\frac{\mathrm{1}}{n^s}}$ in the critical region. In other words, this approach proves that we can place the value s (input) and function  $\zeta \left(s\right)$ in a one-to-one correspondence with the infinite series $1+\frac{1}{2^s}+\frac{1}{3^s}+\textit{{\dots}}$ in the critical region. Moreover, the equality signs in \eqref{2-5-e} prove that the infinite series $1+\frac{1}{2^s}+\frac{1}{3^s}+\textit{{\dots}}$ would necessarily overlap with the analytic continuation of the zeta function.\\

\noindent It is worth mentioning here that the equality signs in \eqref{2-5-e} are concrete proof that the analytic continuation of the Riemann zeta function is the analytic continuation of the Euler product   ($\mathrm\prod_p{{\left(1-{p^{-s}}\right)}^{-1}}$) [\textit{(1), Chap. I}] in the critical region. Thus, when considering \eqref{2-5-e}, it is irrational to assume that establishing a mathematical relation between the Euler product and analytic continuation of the Riemann zeta function in the critical region is wrong because the Euler product does not exist(diverges $\sigmaup>=1$ [\textit{(1), Chap. I}])  in the critical region (i.e., it is irrational to assume that the complex analysis is wrong because $\sqrt{-1}$ does not exist).\\

\noindent Now, for simplicity and clarity, let $\mathrm{\alphaup}=\left(1-\frac{2}{2^s}\right)$ and $\mathrm{\betaup}=\left(1-\frac{2}{2^{1-s}}\right)$. Then, considering \eqref{2-5-e}, we obtain the one-to-one correspondence or bijection between the alternating zeta function $\mathrm{\zetaup }\left(s\right)=\mathrm{\ }\left(\sum^{\mathrm{\infty }}_{n\mathrm{=1}}{\frac{\mathrm{-}{\mathrm{1}}^{n-1}}{n^s}}\right)\ {\left(1-\frac{2}{2^s}\right)}^{-1}$ and the series $\mathop\mathrm{lim}_{b\mathrm{\to }\mathrm{\infty }}\left(\sum^b_{n=1}{\frac{1}{n^s}}\right)$, where $\mathrm{\zetaup }\left(s\right)=\overline{\mathrm{\zetaup }\left(1-s\right)}$.
\begin{align}\label{2-5}
	\begin{split}
		\left(\sum^{\infty }_{n=1}{\frac{-1^{n-1}}{n^s}}\right){\left(\mathrm{\alphaup}\right)}^{-1}=\ &\overline{\left(\sum^{\infty }_{n=1}{\frac{-1^{n-1}}{n^{1-s}}}\right){\left(\mathrm{\betaup}\right)}^{-1}}\\
		=\sum^{\infty }_{n=1}{\left(\frac{1}{n^s}-\frac{2}{{\left(2n\right)}^s}\right){\left(\mathrm{\alphaup}\right)}^{-1}}=\ &\overline{\sum^{\infty }_{n=1}{\left(\frac{1}{n^{1-s}}-\frac{2}{{\left(2n\right)}^{1-s}}\right){\left(\mathrm{\betaup}\right)}^{-1}}}\\
		=\sum^{\infty }_{n=1}{\frac{1}{n^s}\left(\mathrm{\alphaup}\right){\left(\mathrm{\alphaup}\right)}^{-1}}=\ &\overline{\sum^{\infty }_{n=1}{\frac{1}{n^{1-s}}\left(\mathrm{\betaup}\right){\left(\mathrm{\betaup}\right)}^{-1}}}\\ 
		=1+\frac{1}{2^s}+\frac{1}{3^s}{+\dots} =\ &\overline{1+\frac{1}{2^{1-s}}+\frac{1}{3^{1-s}}{+\dots}} \\
		={\mathop\mathrm{lim}_{b\to \infty } \left(\sum^b_{n=1}{\frac{1}{n^s}}\right)\ }=\ &\overline{{\mathop\mathrm{lim}_{b\to \infty } \left(\sum^b_{n=1}{\frac{1}{n^{1-s}}}\right)\ }}\\
	\end{split}
\end{align}

\noindent The main idea of \eqref{2-5} is to prove the one-to-one mapping \eqref{diag1} between the set of finite values $\mathbb{A}=\left\{\zeta \left(s\right) \left|\ s\in \ \mathbb{S}\right.\right\}=\left\{ \overline{\ \zeta \left(1-s\right)}\ \left|\ s\in \mathbb{S}\right.\right\}$ and set of infinite series $\mathbb{B}=\left\{{\mathop\mathrm{lim}_{b\to \infty } \left(\sum^b_{n=1}{\left(\frac{1}{n^s}\right)}\right)\ }\left|\ s\in \ \mathbb{S}\right.\ \right\}=\left\{\ {\mathop\mathrm{lim}_{b\to \infty } \left(\sum^b_{n=1}{\left(\overline{\frac{1}{n^{1-s}}}\right)}\right)\ }\left|\ s\in \ \mathbb{S}\right.\ \right\}$, where $\mathbb{S}=\left\{s\ \left|\ 0<\mathfrak{R}\left(s\right)<1,\ \zeta \left(s\right)=\overline{\zeta \left(1-s\right)}\right.\right\}$.

\begin{equation}\label{diag1}
	\begin{CD}
		\mathbb{B}=\left\{ 1+\frac{1}{2^s}+\frac{1}{3^s}{+\dots}\left|\ s\in \ \mathbb{S}\right.\ \right\}                @=    \left\{ \overline{1+\frac{1}{2^{1-s}}+\frac{1}{3^{1-s}}{+\dots}}\left|\ s\in \ \mathbb{S}\right.\ \right\}  \\
		@VV{\eqref{2-5-e}}V                                @VV{\eqref{2-5-e}}V\\
		\mathbb{A}=\left\{\zeta \left(s\right) \left|\ s\in \ \mathbb{S}\right.\right\}  @=   \left\{ \overline{\ \zeta \left(1-s\right)}\ \left|\ s\in \mathbb{S}\right.\right\}\\
	\end{CD}
\end{equation}
\\
\noindent Considering $ \mathrm{\zetaup }\left(s\right)=\overline{\mathrm{\zetaup }\left(1-s\right)} \mathrm{\Leftrightarrow }\mathrm{\zetaup }\left(s\right)-\overline{\mathrm{\zetaup }\left(1-s\right)}=0$, \eqref{2-5-e} and \eqref{2-5} can be rewritten as
\begin{align}\label{2-5-1}
	\begin{split}
		&\zeta \left(s\right)-\overline{\zeta \left(1-s\right)}\mathrm{=0}\\
		\mathrm{\Leftrightarrow }&\left(\sum^{\infty }_{n=1}{\frac{-1^{n-1}}{n^s}}\right){\left(\alphaup\right)}^{-1}-\overline{\left(\sum^{\infty }_{n=1}{\frac{-1^{n-1}}{n^{1-s}}}\right){\left(\betaup\right)}^{-1}}\mathrm{=0}\\
		\mathrm{\Leftrightarrow }&\sum^{\infty }_{n\mathrm{=1}}{\left(\left({\frac{-1^{n-1}}{n^s}}\right){\left(\alphaup\right)}^{-1}-\overline{\left({\frac{-1^{n-1}}{n^{1-s}}}\right){\left(\betaup\right)}^{-1}}\right)}\mathrm{=0}\\
		\mathrm{\Leftrightarrow }&\sum^{\infty }_{n\mathrm{=1}}{\left(\left(\frac{\mathrm{1}}{n^s}-\frac{2}{{\left(2n\right)}^s}\right){\left(\alphaup\right)}^{-1}-\overline{\left(\frac{\mathrm{1}}{n^{1-s}}-\frac{2}{{\left(2n\right)}^{1-s}}\right){\left(\betaup\right)}^{-1}}\right)}\mathrm{=0}\\
				\mathrm{\Leftrightarrow }& \sum^{\mathrm{\infty }}_{n=1}{\left(\frac{1}{n^s}\alphaup{\left(\alphaup\right)}^{-1}-\overline{\frac{1}{n^{1-s}}\betaup{\left(\betaup\right)}^{-1}}\right)}\mathrm{=0}
	\end{split}
\end{align}

\noindent Consequently, since ${\alphaup}{\left(\alphaup\right)}^{-1}=1=\overline{{\betaup}{\left(\betaup\right)}^{-1}}$ clearly, we have
\begin{align}\label{2-5-2}
	\begin{split}
		\mathrm{\zetaup }\left(s\right)-\overline{\mathrm{\zetaup }\left(1-s\right)}\mathrm{\ }\mathrm{=0=} {\mathop{\mathrm{lim}}_{b\mathrm{\to }\mathrm{\infty }} \left(\sum^b_{n=1}{\left(\frac{1}{n^s}-\left(\overline{\frac{1}{n^{1-s}}}\right)\right)}\ \right)\ , \ s\in \ \mathbb{S} }
	\end{split}
\end{align}

\noindent In the following part, we use the fact that if the analytic continuation of $\mathrm{f}\left(s_1\right)$ is equal to the analytic continuation of $\mathrm{f}\left(s_2\right)$, then $\mathrm{f}\left(s_1\right)-\mathrm{f}\left(s_2\right)\mathrm{=0}$ (and vice versa) to prove \eqref{2-5-2} by an alternative method. Then, we use \eqref{2-5-2} to transform \eqref{2-4-3} into the desired exponential equality.\\

\noindent Suppose that ${\mathop{\mathrm{lim}}_{{b}\mathrm{\to }\mathrm{\infty }} \left(\sum^{{b}}_{{n}=1}{\left(\frac{1}{{{n}}^{{s}}}-\left(\overline{\frac{1}{{{n}}^{1-\mathrm{s}}}}\right)\right)}\right)\ }\mathrm{=0}$. Therefore, considering \eqref{2-5-e} and \eqref{2-5}, we obtain $\mathrm{\zetaup }\left(s\right)\mathrm{-}\overline{\mathrm{\zetaup }\left(\mathrm{1-}s\right)}\mathrm{=0}$. Then, because $\mathrm{\zetaup }\left(s\right)-\overline{\mathrm{\zetaup }\left(1-s\right)}\mathrm{=0}\mathrm{\Leftrightarrow }\mathrm{\zetaup }\left(s\right)=\overline{\mathrm{\zetaup }\left(1-s\right)}\mathrm{\ }$, we can readily deduce that
\begin{align}\label{2-5-3}
	\begin{split}
		\mathrm{\zetaup }\left(s\right)=\overline{\mathrm{\zetaup }\left(1-s\right)}\mathrm{\Leftrightarrow }{\mathop{\mathrm{lim}}_{b\mathrm{\to }\mathrm{\infty }} \left(\sum^b_{n=1}{\left(\frac{1}{n^s}-\overline{\left(\frac{1}{n^{1-s}}\right)}\right)}\right)\mathrm{=0}}
	\end{split}
\end{align}
\noindent Finally, by substituting ${\mathop{\mathrm{lim}}_{{b}\mathrm{\to }\mathrm{\infty }} \left(\sum^{{b}}_{{n}=1}{\left(\frac{1}{{{n}}^{{s}}}-\left(\overline{\frac{1}{{{n}}^{1-\mathrm{s}}}}\right)\right)}\right)\ }\mathrm{=0}$ into \eqref{2-4-3}, we obtain $0={\mathop{\mathrm{lim}}_{b\mathrm{\to }\mathrm{\infty }} \left(\frac{b^{1-s}}{1-s}-\overline{\left(\frac{b^s}{s}\right)}\right)\ }$ for the zeros of the zeta function. Thus, we have
\begin{align}\label{2-6}
	\begin{split}
		\mathrm{\zetaup }\left({s}\right)=0\ \ {\Rightarrow }\mathop{\mathrm{lim}}_{{b}\mathrm{\to }\mathrm{\infty }}\left(\frac{{{b}}^{1-{s}}}{1-s}-\frac{{{b}}^{\overline{s}}}{\overline{{s}}}\right)=0\mathrm{\ }\mathrm{\Rightarrow }{\mathop{\mathrm{lim}}_{{b}\mathrm{\to }\mathrm{\infty }} \left(\frac{{{b}}^{\mathrm{1-}\mathrm{\sigmaup }}}{\left|{1-s}\right|}-\frac{{{b}}^{\mathrm{\sigmaup }}}{\left|{s}\right|}\right)\ }=0
	\end{split}
\end{align}

\noindent Now consider a simple exponential equality of the form ${\mathop{\mathrm{lim}}_{x\mathrm{\to }\mathrm{\infty }} \frac{x^{\mathrm{\alphaup }}}{\mathrm{\alphaup }}\ }-{\mathop{\mathrm{lim}}_{x\mathrm{\to }\mathrm{\infty }} \frac{x^{\mathrm{\betaup }}}{\mathrm{\betaup }}\ }=0$. We easily see that this equality is satisfied if and only if $\mathrm{\alphaup }\mathrm{=}\mathrm{\betaup }$. Therefore, it is easy to observe that \eqref{2-6} is satisfied if and only if $\mathrm{\sigmaup }\mathrm{=}\frac{\mathrm{1}}{\mathrm{2}}$. We can readily assert this with the following methods: 

\noindent Consider the chronological progress tables below$\ (0<\mathrm{\sigmaup }\mathrm{\ }\mathrm{<1<}\left|t\right|):$ \\

\noindent If $\mathrm{1-}\mathrm{\sigmaup }>\ \mathrm{\sigmaup }\ \left(\mathrm{\Leftrightarrow }\ \mathrm{\ }\mathrm{\sigmaup }\mathrm{\ }\mathrm{<}\mathrm{\ }\frac{\mathrm{1}}{\mathrm{2}}\right)\ $, we have 
\begin{center}
	\begin{tabular}{|p{0.7in}|p{0.7in}|p{0.7in}|p{0.7in}|p{1.2in}|} \hline 
		${b}$\textit{} & $2$\textit{} & $3$\textit{} & \textit{{\dots}} & $\mathrm{lim}_{b\mathrm{\to }\mathrm{\infty }}b\ $\textit{} \\ \hline 
		$\frac{b^{\mathrm{1-}\mathrm{\sigmaup }}}{\left|\mathrm{1-}s\right|}-\frac{b^{\mathrm{\sigmaup }}}{\left|s\right|}$\textit{} & $\frac{2^{\mathrm{1-}\mathrm{\sigmaup }}}{\left|\mathrm{1-}s\right|}-\frac{2^{\mathrm{\sigmaup }}}{\left|s\right|}$\textit{} & $\frac{3^{\mathrm{1-}\mathrm{\sigmaup }}}{\left|\mathrm{1-}s\right|}-\frac{3^{\mathrm{\sigmaup }}}{\left|s\right|}$\textit{} & \textit{{\dots}} & ${\mathop{\mathrm{lim}}_{b\mathrm{\to }\mathrm{\infty }} \left(\frac{b^{\mathrm{1-}\mathrm{\sigmaup }}}{\left|\mathrm{1-}s\right|}-\frac{b^{\mathrm{\sigmaup }}}{\left|s\right|}\right)\ }$\textit{} \\ \hline 
		\textit{} & $\boldsymbol{>}$\textbf{\textit{}} & $\boldsymbol{>}$\textbf{\textit{}} & $\boldsymbol{>}$\textbf{\textit{}} & $\boldsymbol{>}$\textbf{\textit{}} \\ \hline 
		0\textit{} & 0\textit{} & 0\textit{} & 0\textit{} & $0$ \\ \hline 
	\end{tabular}
\end{center}

\noindent Note that $\frac{b^{\mathrm{1-}\mathrm{\sigmaup }}}{\left|\mathrm{1-}s\right|}-\frac{b^{\mathrm{\sigmaup }}}{\left|s\right|}$ cannot be equal to zero as $b$ increases. Proceeding in a similar fashion as the previous proof, we observe that ${\mathop{\mathrm{lim}}_{b\mathrm{\to }\mathrm{\infty }} \left(\frac{b^{\mathrm{1-}\mathrm{\sigmaup }}}{\left|\mathrm{1-}s\right|}-\frac{b^{\mathrm{\sigmaup }}}{\left|s\right|}\right)\ }$ cannot be equal to zero if $\mathrm{1-}\mathrm{\sigmaup }\ \mathrm{<}\mathrm{\ }\mathrm{\sigmaup }\mathrm{\ }\left(\mathrm{\Leftrightarrow }\mathrm{\ }\mathrm{\sigmaup }\mathrm{\ }>\ \frac{\mathrm{1}}{\mathrm{2}}\right)$.\\

\noindent If $\mathrm{1-}\mathrm{\sigmaup }\ \mathrm{=}\mathrm{\ }\mathrm{\sigmaup }\mathrm{\ }$, we have 
\begin{center}
	\begin{tabular}{|p{0.7in}|p{0.7in}|p{0.7in}|p{0.7in}|p{1.2in}|} \hline 
		${b}$\textit{} & $2$\textit{} & $3$\textit{} & \textit{{\dots}} & $\mathrm{lim}_{b\mathrm{\to }\mathrm{\infty }}b\ $\textit{} \\ \hline 
		$\frac{b^{\mathrm{1-}\mathrm{\sigmaup }}}{\left|\mathrm{1-}s\right|}-\frac{b^{\mathrm{\sigmaup }}}{\left|s\right|}$\textit{} & $\frac{2^{\mathrm{1-}\mathrm{\sigmaup }}}{\left|\mathrm{1-}s\right|}-\frac{2^{\mathrm{\sigmaup }}}{\left|s\right|}$\textit{} & $\frac{3^{\mathrm{1-}\mathrm{\sigmaup }}}{\left|\mathrm{1-}s\right|}-\frac{3^{\mathrm{\sigmaup }}}{\left|s\right|}$\textit{} & \textit{{\dots}} & ${\mathop{\mathrm{lim}}_{b\mathrm{\to }\mathrm{\infty }} \left(\frac{b^{\mathrm{1-}\mathrm{\sigmaup }}}{\left|\mathrm{1-}s\right|}-\frac{b^{\mathrm{\sigmaup }}}{\left|s\right|}\right)\ }$\textit{} \\ \hline 
		\textit{} & $\boldsymbol{=}$\textbf{\textit{}} & $\boldsymbol{=}$\textbf{\textit{}} & $\boldsymbol{=}$\textbf{\textit{}} & $\boldsymbol{=}$\textbf{\textit{}} \\ \hline 
		0\textit{} & 0 & 0 & 0\textit{} & $0$\textit{} \\ \hline 
	\end{tabular}
\end{center}

\noindent Note that $\frac{b^{\mathrm{1-}\mathrm{\sigmaup }}}{\left|\mathrm{1-}s\right|}-\frac{b^{\mathrm{\sigmaup }}}{\left|s\right|}$ is equal to zero for all values of b $\mathrm{>}$ 1 $\left(b\mathrm\ {\in }\mathrm{\ }\mathbb{N}\right)$. Thus, we can conclude that ${\mathop{\mathrm{lim}}_{b\mathrm{\to }\mathrm{\infty }} \frac{b^{\mathrm{1-}\mathrm{\sigmaup }}}{\left|\mathrm{1-}s\right|}\ }\mathrm{-}{\mathop{\mathrm{lim}}_{b\mathrm{\to }\mathrm{\infty }} \frac{b^{\mathrm{\sigmaup }}}{\left|s\right|}\ }=0$ if and only if $\mathrm{1-}\mathrm{\sigmaup }\ \mathrm{=}\mathrm{\ }\mathrm{\sigmaup }\left(\mathrm{\Leftrightarrow }\ \mathrm{\sigmaup }=\frac{1}{2}\right).$\\

\noindent Hence, by considering \eqref{2-6}, we have
\begin{align*}
	\begin{split}
		\mathrm{1-}\mathrm{\ }\mathrm{\sigmaup }\ \boldsymbol{\mathrm{>}}\mathrm{\ }\mathrm{\sigmaup }\ &\mathrm{\Rightarrow }{\mathop{\mathrm{lim}}_{b\mathrm{\to }\mathrm{\infty }} \left(\frac{b^{\mathrm{1-}\mathrm{\sigmaup }}}{\left|\mathrm{1-}s\right|}-\frac{b^{\mathrm{\sigmaup }}}{\left|s\right|}\right)\ }\mathrm{\ }\mathrm{=+}\mathrm{\infty }\ \boldsymbol{\mathrm{>}}\ \mathrm{0}\\
		\mathrm{\ }\mathrm{1-}\mathrm{\ }\mathrm{\sigmaup }\ \boldsymbol{\mathrm{<}}\mathrm{\ }\mathrm{\sigmaup }\ &\mathrm{\Rightarrow }{\mathop{\mathrm{lim}}_{b\mathrm{\to }\mathrm{\infty }} \left(\frac{b^{\mathrm{1-}\mathrm{\sigmaup }}}{\left|\mathrm{1-}s\right|}-\frac{b^{\mathrm{\sigmaup }}}{\left|s\right|}\right)\ }\mathrm{=-}\mathrm{\infty }\ \boldsymbol{\mathrm{<}}\ \mathrm{0}
	\end{split}
\end{align*}

\noindent Thus, we conclude that if $\mathrm{\sigmaup }<\frac{1}{2}\ $ or $\ \mathrm{\sigmaup }>\frac{1}{2}$, $\ $the equality ${\mathop{\mathrm{lim}}_{b\mathrm{\to }\mathrm{\infty }} \left(\frac{b^{\mathrm{1-}\mathrm{\sigmaup }}}{\left|\mathrm{1-}s\right|}-\frac{b^{\mathrm{\sigmaup }}}{\left|s\right|}\right)\ }=0$ cannot be satisfied.\\

\noindent We observe that \eqref{2-6} holds if and only if $\mathrm{\sigmaup }=\frac{1}{2}$. In other words, $\mathrm{\sigmaup }\mathrm{\neq }\frac{1}{2}\mathrm{\Rightarrow }\mathrm{\zetaup }\left(s\right)\mathrm{\neq }\overline{\mathrm{\zetaup }\left(1-s\right)}$. Hence, we conclude that
\begin{align}
	\begin{split}\label{6}
		\mathrm{\sigmaup\neq\frac{1}{2}\Leftrightarrow\zetaup\left(s\right)\neq0}
	\end{split}
\end{align}

\noindent This completes the proof that the real part of the nontrivial zeros of the zeta function is equal to 1/2. \\

\noindent Over the years, hundreds of mathematical theories have been built upon the assumption that Riemann's last theorem is true. Therefore, considerable efforts have been made by several of the best mathematical minds around the world to protect the legitimacy of these theories. However, in this work we have finally proved this famous theorem that had resisted all efforts to be proven for over one and a half centuries.

\begin{figure}[h]
	\begin{flushright}
							\includegraphics[width=.1\textwidth]{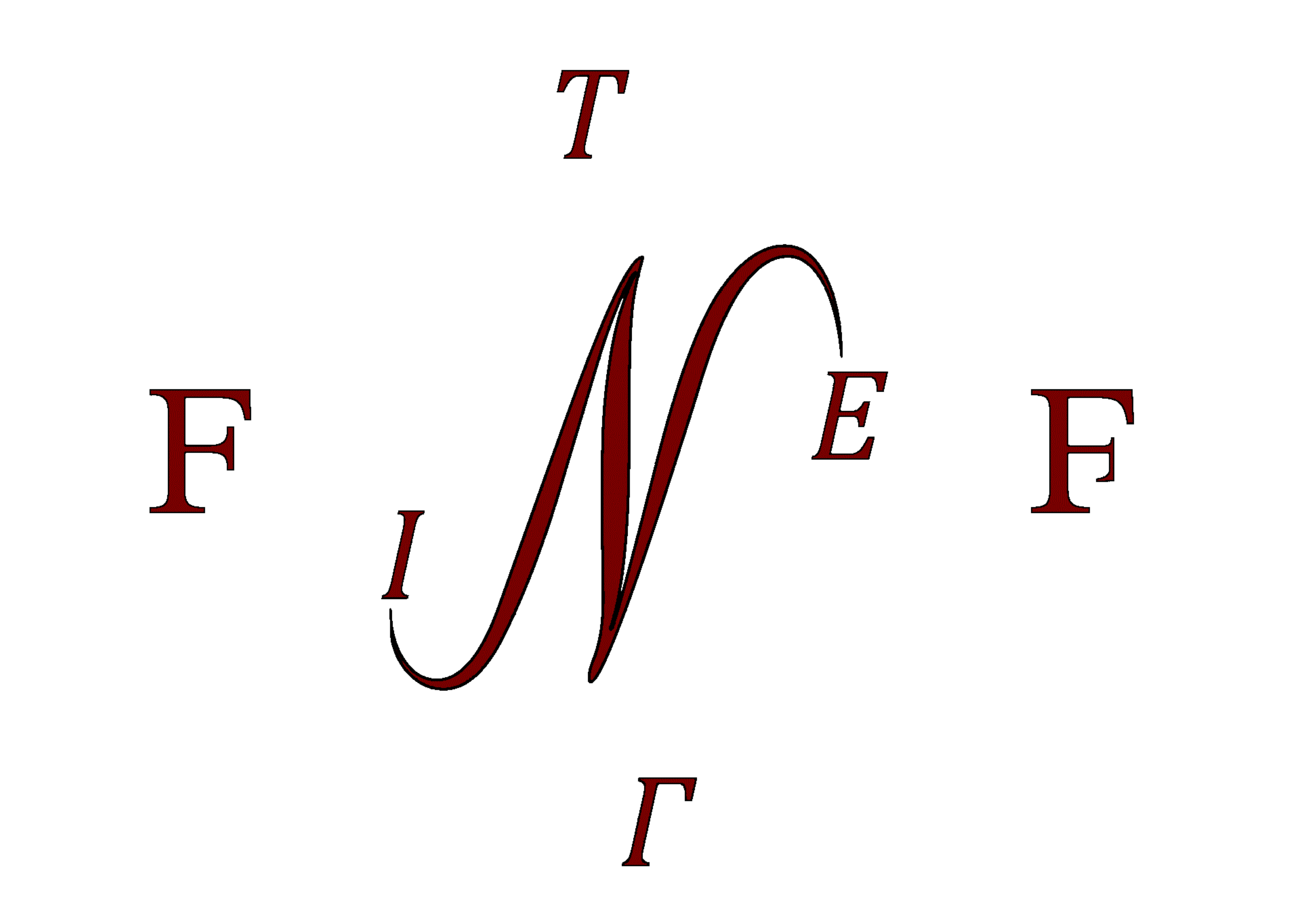}
	\end{flushright}
\end{figure}
%\pagebreak
\section*{Declarations}
\section*{Funding} 
Not applicable.

\section*{Conflict of interest}
The author declares that he has no affiliations with or involvement in any organization or entity with any financial interest (such as honoraria, educational grants, participation in speakers’ bureaus, membership, employment, consultancies, stock ownership, or other equity interest, and expert testimony or patent-licensing arrangements), or non-financial interest (such as personal or professional relationships, affiliations, knowledge or beliefs) in the subject matter or materials discussed in this manuscript.

\section*{Availability of data and material}
All the data and materials are available upon request.

\section*{Code availability}
All the code are available upon request.

\section*{Author contributions}
Not applicable.


\begin{thebibliography}{}



\bibitem{Ref1}
\noindent  \textit{Titchmarsh, E.C., Wolffe, A.: The Theory of the Riemann Zeta-Function. Clarendon Press, United Kingdom (1986).}\\

\bibitem{Ref2}
\noindent  \textit{Edwards, H.M.: Riemann's Zeta Function. Dover Publications, United States (2001).}\\

\bibitem{Ref3}
\noindent \textit{Kadry, S.: Mathematical Formulas for Industrial and Mechanical Engineering. Elsevier Science, India (2014).}\\

\bibitem{Ref4}
\noindent  \textit{Dimitris, K.: The Distribution of Prime Numbers. American Mathematical Society, United States (2020).}\\ 

\bibitem{Ref5}
\noindent  \textit{Cochrane, L.: Adelard of Bath: The First English Scientist. British Museum Press, London (1994).}\\ 

\bibitem{Ref6}
\noindent  \textit{Berndt, B.C.: Ramanujan’s Notebooks: Part I. Springer, New York (2012).}\\ 

\bibitem{Ref7}
\noindent  \textit{Berndt, B.C.: Ramanujan’s Notebooks: Part II. Springer, New York (2012).}\\

\bibitem{Ref8}
\noindent  \textit{Berndt, B.C.: Ramanujan’s Notebooks: Part III. Springer, New York (2012).}\\

\bibitem{Ref9}
\noindent  \textit{Berndt, B.C.: Ramanujan’s Notebooks: Part IV. Springer, New York (2012).}\\

\bibitem{Ref11}
\noindent \textit{Berndt, B.C.: Ramanujan’s Notebooks: Part V. Springer, New York (2012).}\\

\bibitem{Ref10}
\noindent \textit{Choi, J., Srivastava, H.M.: Zeta and Q-Zeta Functions and Associated Series and Integrals. Elsevier Science, Netherlands (2012). }\\

\bibitem{Ref12}
\noindent  \textit{Ghusayni, B.: Classical Real Analysis. Badih Ghusayni, N.p. (2009).}\\

\bibitem{Ref13}
\noindent  \textit{Ghusayni, B.: Number Theory from an Analytic Point of View. Badih Ghusayni, N.p. (2004).}\\

\bibitem{Ref14}
\noindent  \textit{Hardy, G.H., Wright, E.M., Wiles, A., Silverman, J., Heath-Brown, R.: An Introduction to the Theory of Numbers. Oxford University Press, Oxford (2008).}\\

\bibitem{Ref15}
\noindent  \textit{Ivic, A.: The Riemann Zeta-Function: Theory and Applications. Dover Publications, United States (2012).}\\

\bibitem{Ref16}
\noindent  \textit{Polchinski, J.: String Theory: Volume 1, An Introduction to the Bosonic String. Cambridge University Press, Cambridge (1998).}\\

\bibitem{Ref17}
\noindent  \textit{Voronin, S.M., Koblitz, N., Karatsuba, A.A.: The Riemann Zeta-Function. De Gruyter, Germany (2011).}\\


\end{thebibliography}
\end{document}